# Real-time identification of the current density profile in the JET Tokamak: method and validation

D. Mazon, J. Blum, C. Boulbe, B. Faugeras, A. Boboc, M. Brix, P. De Vries, S. Sharapov, L. Zabeo

*Abstract*—The real-time reconstruction of the plasma magnetic equilibrium in a Tokamak is a key point to access high performance regimes. Indeed, the shape of the plasma current density profile is a direct output of the reconstruction and has a leading effect for reaching a steady-state high performance regime of operation. In this paper we present the methodology followed to identify numerically the plasma current density in a Tokamak and its equilibrium. In order to meet the real-time requirements a C++ software has been developed using the combination of a finite element method, a nonlinear fixed point algorithm associated to a least square optimization procedure. The experimental measurements that enable the identification are the magnetics on the vacuum vessel, the interferometric and polarimetric measurements on several chords and the motional Stark effect. Details are given about the validation of the reconstruction on the JET tokamak, either by comparison with 'off-line' equilibrium codes or real time software computing global quantities.

## I. Introduction

In present days tokamaks the shape of the plasma boundary is routinely identifiable in real-time in less than few milliseconds using a set of magnetic and diamagnetic coils spread around the vessel [1]. This information is mainly used for controlling the plasma shape in real-time during a plasma discharge using coils current in a feedback control loop. The idea is to achieve a required shape and to maintain it in a stationary manner in order to avoid for example sudden termination of the plasma when the plasma touches the first wall. In JET the so-called XLOC code is used routinely for plasma shape control. Based on this JET flux boundary code confinement parameters are deduced like the diamagnetic energy, the internal inductance and plasma separatrix geometry in less than 1ms. But with this algorithm it is not possible to compute the internal magnetic flux configuration which is needed if we want to analyze the phenomenon occurring in the interior of the plasma. In this case the only way to get access to the current density profile is to use off-line codes that can compute accurately the profile but with no possibility to act on it. This is rather a strong limitation because we know from the analysis performed that the shape of the current density profile is one of the key element to enhanced the plasma performances. We have seen in particular that non monotonic current density profiles can trigger enhanced particles and heat confinement [2]. On top of this the current density profile presents a resistive diffusion time and react with a delay on any variation of the current drive systems. So it is clear that by controlling in real time such a profile we ensure stability but also performances [3, 4]. The problem of the equilibrium of a plasma in a Tokamak is a free boundary problem in which the plasma boundary is defined as the last closed magnetic flux surface. Inside the plasma, the equilibrium equation in an axisymmetric configuration is called the Grad-Shafranov equation [5, 6]. The right hand side of this equation is a non-linear source which represents the toroidal component of the plasma current density. The goal of a real-time equilibrium code is to identify not only the plasma boundary but also the flux surface geometry outside and inside of the plasma, the current density profile and derive the safety factor 'q' and other important parameters from the obtained equilibrium. In order to meet the real-time requirements, a new version of the code called Equinox has been designed and implemented in C++. The code relies on tokamak specific software providing flux values on the first wall of the vacuum vessel. By means of least-square minimization of the difference between magnetic measurements and the simulated ones the code identifies the source term of the non linear Grad-Shafranov equation. The finite element solver uses triangles mesh, the calculation being limited to the vacuum chamber. A careful implementation leads to execution time less than 60ms per iteration on a 2GHz PC, complemented with excellent robustness and very good precision (+/- 1cm) of plasma boundary for an equilibrium code. In the following, the next section is devoted to the mathematical modeling of the equilibrium problem in axisymmetric configurations. Then

D. Mazon is with the Association EURATOM-CEA, CEA Cadarache, DSM/IRFM 13108 St Paul Lez durance Cedex France phone: 00 33 4 42 25 48 53; fax: 00 33 4 42 25 26 61; (e-mail: dmazon@jet.uk).

C. Boulbe is with the Laboratoire J-A Dieudonné (UMR 66 21), Université de Nice Sophia-Antipolis, CNRS Parc Valrose 06108 Nice Cedex 02 France (e-mail: boulbe@unice.fr).

B. Faugeras is with the Laboratoire J-A Dieudonné (UMR 66 21), Université de Nice Sophia-Antipolis, CNRS Parc Valrose 06108 Nice Cedex 02 France (e-mail: faugeras@unice.fr).

J. Blum is with the Laboratoire J-A Dieudonné (UMR 66 21), Université de Nice Sophia-Antipolis, CNRS Parc Valrose 06108 Nice Cedex 02 France (e-mail: jblum@unice.fr).

A. Boboc is with the association EURATOM/UKAEA, Culham science centre Abingdon Oxon OX14 3DB UK. (e-mail: aboboc@jet.uk).

M. Brix is with the association EURATOM/UKAEA, Culham science centre Abingdon Oxon OX14 3DB UK. (e-mail: mbrix@jet.uk).

P. De. Vries is with the association EURATOM/UKAEA, Culham science centre Abingdon Oxon OX14 3DB UK. (e-mail: pvries@jet.uk).

S. Sharapov is with the association EURATOM/UKAEA, Culham science centre Abingdon Oxon OX14 3DB UK. (e-mail: sershar@jet.uk).

L. Zabeo is with the association EURATOM/UKAEA, Culham science centre Abingdon Oxon OX14 3DB UK. (e-mail: lzabeo@jet.uk).

the inverse problem will be addressed. Finally, the validation of the results using a database of 130 discharges with a large variety of magnetic configurations, plasma current and toroidal field strength will be given in the final section. Distinction will be made between the Equinox M version which is using only magnetic measurements as constraints and the Equinox J version which is constrained using internal measurement like interferometry, polarimetry, and Motional Stark effect data. Validation will be presented for the Equinox M version.

## II. THE GRAD-SHAFRANOV EQUATION

### A. Mathematical modeling

In the presence of a magnetic field $B$ the equations governing the plasma equilibrium are first the magnetostatic Maxwell's equations which are satisfied in the whole of space (including the plasma itself):

$$\begin{cases} \nabla \cdot B = 0 \\ \nabla \times (\frac{B}{\mu}) = j \end{cases} \quad (1)$$

where $\mu$ represents the magnetic permeability and $j$ is the current density and second the equilibrium for the plasma itself .which can be written as follow at the resistive time scale [7]:

$$\nabla p = j \times B \quad (2)$$

It is clear from equation (2) that the plasma is in equilibrium when the force $\nabla p$ due to the kinetic pressure is equal to the Lorentz force of the magnetic pressure $j \times B$. As a consequence the field lines and the current lines for a plasma in equilibrium lie on isobaric surfaces. These surfaces, generated by the field lines are called magnetic surfaces. As they need to remain within a bounded volume of space it is clear that they need to have a toroidal topology. These surfaces constitute a set of nested tori. The torus degenerates progressively into a curve which is called the magnetic axis (innermost torus). In a cylindrical coordinate system $(r, z, \phi)$ where $r=0$ is the major axis of the torus, the hypothesis of axial symmetry consists in assuming that the magnetic field $B$ is independent of the toroidal angle $\phi$. The magnetic field is usually decomposed as $B=B_p+B_\phi$ where $B_p = (Br, Bz)$ is the poloidal component and $B_\phi$ is the toroidal component. From equation (1) we can define the poloidal flux $\psi(r,z)$ such that:

$$\begin{cases} B_r = -\frac{1}{r}\frac{\partial \psi}{\partial z} \\ B_z = \frac{1}{r}\frac{\partial \psi}{\partial r} \end{cases} \quad (3)$$

If we note $e_\phi$ the unit vector in the toroidal direction and $f$ the diamagnetic function the poloidal and toroidal component of the magnetic field can been written respectively as in the following formula:

$$\begin{cases} B_p = \frac{1}{r}[\nabla \psi \times e_\phi] \\ B_\phi = \frac{f}{r} e_\phi \end{cases} \quad (4)$$

From equation (4) and the second relation of (1) we obtain the following expression for $j_p$ and $j_\phi$ respectively the poloidal and toroidal component of j:

$$\begin{cases} j_p = \frac{1}{r}\left[\nabla\left(\frac{f}{\mu}\right) \times e_\phi\right] \\ j_\phi = (-\Delta^*\psi)e_\phi \end{cases} \quad (5)$$

The linear elliptic operator $\Delta^*$ being defined by

$$\Delta^* = \frac{\partial}{\partial r}\left(\frac{1}{\mu r}\frac{\partial}{\partial r}\right) + \frac{\partial}{\partial z}\left(\frac{1}{\mu r}\frac{\partial}{\partial z}\right) \quad (6)$$

In the plasma region relation (2) implies that $\nabla p$ and $\nabla \psi$ are collinear and therefore $p$ is constant on each magnetic surface. This can be put in the following mathematical form:

$$p = p(\psi) \quad (7)$$

In the same way combining the expression (5) and (2) implies immediately that $\nabla p$ and $\nabla f$ are also collinear and therefore f is likewise constant on each magnetic surface

$$f = f(\psi) \quad (8)$$

The equilibrium relation (2) combined with the expression (4) and (5) for B and j implies that:

$$\nabla p = -\frac{\Delta^*\psi}{r}\nabla\psi - \frac{f}{\mu_0 r^2}\nabla f \quad (9)$$

which leads to the so-called Grad-Shafranov equation:

$$-\Delta^*\psi = rp'(\psi) + \frac{1}{\mu_0 r}(ff')(\psi) \quad (10)$$

with $\Delta^*$ the linear operator defined by (6) in which $\mu$ is equal to the magnetic permeability $\mu_0$ of the vacuum. From (5) it is clear that the right hand side of the equation (10) represents the toroidal component of the plasma current density. It involves functions $p(\psi)$ and $f(\psi)$ which are not directly measured inside the plasma. We can also note that in the vacuum where no plasma current is present the magnetic flux $\psi$ satisfies:

$$-\Delta^*\psi = 0 \quad (11)$$

The equilibrium of a plasma in a domain $\Omega$ representing the vacuum region is a free boundary problem. The plasma free boundary is defined at JET as being a magnetic separatrix (hyperbolic line with an X-point X) the region $\Omega p$ containing the plasma is defined as

$$\Omega_p = \{x \in \Omega, \psi(x) \geq \psi\} \quad (12)$$

where $\psi_b = \psi(X)$ in the X point configuration (see Fig.1). Assuming Dirichlet boundary conditions, h, are given on $\Gamma = \partial\Omega$ which is the poloidal cross section of the vacuum vessel, the final equations governing the behaviour of $\psi(r,z)$ inside the vacuum vessel are:

$$\begin{cases} -\Delta^* \psi = \left[ \dfrac{r}{R_0} A(\overline{\psi}) + \dfrac{R_0}{r} B(\overline{\psi}) \right] \chi\Omega_p & \text{in } \Omega \\ \psi = h & \text{on } \Gamma \end{cases} \quad (13)$$

with

$$A(\overline{\psi}) = R_0 p'(\overline{\psi}) \text{ and } B(\overline{\psi}) = \dfrac{1}{\mu_0 R_0}(ff')(\overline{\psi}) \quad (14)$$

where the normalized flux is introduced so that A and B are defined on the interval [0,1]:

$$\overline{\psi} = \dfrac{\psi - \max\limits_{\Omega_p} \psi}{\psi_b - \max\limits_{\Omega_p} \psi} \quad (15)$$

$\chi\Omega_p$ is the characteristic function of $\Omega_p$.

### B. Statement of the inverse problem

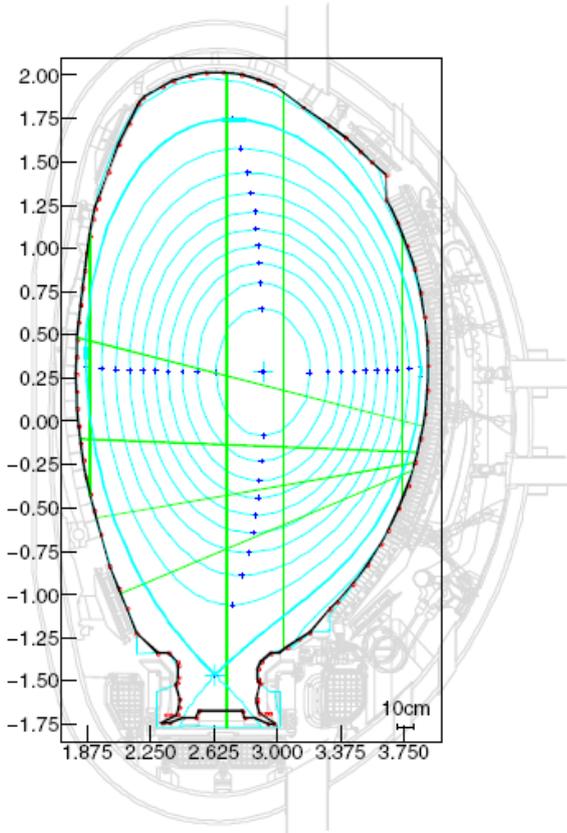

Fig. 1. JET Vessel and plasma boundary (blue thick line). Vertical and horizontal lines reprensent the polarimetry lines of sight.

In order to find the plasma equilibrium we need to solve the non linear two-dimensional partial differential equation (10). The right hand side of this equation is composed of two functions representing the pressure function p and the diamagnetic function f. The numerical identification problem is formulated as a least-square minimization based on available measurements with a Tikhonov regularization [9]. The experimental measurements that enable the identification are the magnetics on the vacuum vessel, the interferometric and polarimetric measures on several chords and the motional Stark effect. For the magnetic measurements the flux loops give the poloidal flux on particular nodes $M_i$ such that $\psi(M_i)=h_i$ on $\Gamma$

$$\psi(M_i) = h_i \text{ on } \Gamma \quad (17)$$

Thanks to an interpolation (performed by XLOC at JET) between the points $M_i$ these measurements provide the Dirichlet boundary condition h. The magnetic probes give the component of the magnetic poloidal field which is tangent to the vacuum vessel

$$\dfrac{1}{r}\dfrac{\partial \psi}{\partial n}(N_i) = g_i \quad (18)$$

The interferometric measurements give the density integrals over the chords $C_i$

$$\int_{C_i} n_e(\overline{\psi}) dl = \beta_i \quad (19)$$

$n_e$ represents the electronic density which is approximately constant on each flux line. The polarimetric measurements give the Faraday rotation of the angle of infrared radiation crossing the section of the plasma along the same chords $C_i$:

$$\int_{C_i} n_e(\overline{\psi}) B_{//} dl = \int_{C_i} \dfrac{n_e(\overline{\psi})}{r} \dfrac{\partial \psi}{\partial n} dl = \alpha_i \quad (20)$$

the component of the poloidal field tangent to $C_i$ is $B_{//}$ and d/dn represents the normal derivative of $\psi$ with respect to $C_i$. The Motional Stark effect (MSE) angle $\gamma_i$ is taken at different points $x_i=(r_i,z_i)$:

$$\tan(\gamma_i) = \dfrac{a_1 B_r + a_2 B_z + a_3 B_\phi}{a_4 B_r + a_5 B_z + a_6 B_\phi} \quad (21)$$

The problem is thus resumed to find a solution that minimizes the cost function defined as:

$$J(A,B,n_e) = J_0 + K_1 J_1 + K_2 J_2 + K_3 J_3 + J_\varepsilon \quad (22)$$

with

$$J_0 = \sum_i \left( \dfrac{1}{r}\dfrac{\partial \psi}{\partial n}(N_i) - g_i \right)^2$$

$$J_1 = \sum_i \left( \int_{C_i} \dfrac{n_e}{r} \dfrac{\partial \psi}{\partial n} dl - \alpha_i \right)^2 \quad (23)$$

$$J_2 = \sum_i \left( \int_{C_i} n_e dl - \beta_i \right)^2$$

$$J_3 = \sum_i (mse_i - \gamma_i)^2$$

with mse the reconstructed measurement and $K_1$ to $K_3$ the weighting parameters enabling to give more or less importance to the corresponding experimental measurements [8]. The inverse problem of the determination of A and B is ill-posed. Hence a regularization procedure has to transform it into a well-posed one [9]. The Tikhonov regularization term $J_\varepsilon$ constrains the function A, B and $n_e$ to be smooth enough and reads:

$$J_\varepsilon = \varepsilon_1 \int_0^1 [A''(x)]^2 dx + \varepsilon_2 \int_0^1 [B''(x)]^2 dx + \varepsilon_3 \int_0^1 [n_e''(x)]^2 dx \quad (24)$$

where $\varepsilon_1$, $\varepsilon_2$ and $\varepsilon_3$ are the regularizing parameters. It is worth mentioning that the electronic density $n_e$ does not intervene in equation (13). However as soon as we want to use the polarimetric measurements it is necessary to include $n_e$ in the identification procedure.

*C. Identification and algorithm*

Equation (13) is solved using a finite element method [10]. The following formulation of this equation forms the basis of the finite element method:

Find $\psi$ such that $\psi = h$ on $\Gamma$ and

$$\int_\Omega \frac{1}{\mu_0 r} \nabla \psi \cdot \nabla v \, dx = \int_\Omega \left[ \frac{r}{R_0} A(\overline{\psi}) + \frac{R_0}{r} B(\overline{\psi}) \right] v \, dx \quad (25)$$

The unknown functions $A$, $B$, $n_e$ are approximated by decomposition in a reduced basis. The basis can be made of different types of functions (polynomials, B-splines, wavelets etc) [11]. In our case we choose B-splines. Let $u$ be the vector which contains the coordinates of $A$, $B$ and $n_e$ in the chosen basis. The Picard type (fixed point) algorithm is then used to solve iteratively the inverse and direct problem. The discretisation of the equation (25) can be written as:

$$\tilde{K} \psi = D(\psi) u + h \quad (26)$$

where $D$ is the plasma current matrix, $\tilde{K}$ is the stiffness matrix and $h$ is due to the Dirichlet boundary conditions. The discrete inverse optimization problem is to find $u$ minimizing the cost function which can be written as

$$J(u) = \|C(\psi)\psi - k\|^2 + u^T \Lambda u \quad (27)$$

while $\psi$ satisfies (26). The quantity $C(\psi)\psi$ represents the outputs of the model, $k$ the experimental measurements, $C(\psi)$ is the observation operator. The matrix $\Lambda$ represents the regularization terms. In order to solve this problem we use an iterative algorithm based on fixed point iterations. At the $n^{th}$ iteration $\psi_n$ and $u_n$ are given. The non linear mapping between $\psi(u)$ and $u$ is given by the relation:

$$\psi = \tilde{K}^{-1}[D(\psi_n) u + h] \quad (28)$$

and the cost function to be minimized is given by

$$J(u) = \|C(\psi_n)\psi - k\|^2 + u^T \Lambda u \quad (29)$$

This last equation is used to determine $u_{n+1}$. Then fixed point iterations for equation (28) enable to find $\psi_{n+1}$.

$$\psi_{n+1} = \tilde{K}^{-1}[D(\psi_n) u_{n+1} + h] \quad (30)$$

Since the algorithm is initialized from the equilibrium at the previous time step two or three iterations are usually enough to ensure convergence. This leads to a very efficient algorithm. In the real-time version of the code the number of fixed point iteration is set to two, taking into account that the equilibrium does not change so much between two consecutive time instances.

Validation of the equilibrium reconstruction

The code installed on JET relies on the boundary code reconstruction XLOC [1] which provides the total plasma current, the toroidal magnetic field and the magnetic flux values and poloidal magnetic field on the first wall of the vacuum vessel. This improves the portability of the code to other machines since we are not asking for the plasma boundary itself. From a practical point of view the boundary codes are traditionally required to give an accurate plasma boundary even if we can observe local oscillations due to high degree polynomial extrapolation used in those codes. Hopefully, while using the boundary conditions on the first wall we use the values where they are the most accurate, while hiding tokamak magnetic measurement specific issues in boundary codes. We can see in Fig. 2 the kind of grid that we choose and the nodes at the boundary where the inputs are computed from the XLOC code. The validation of the Equinox code has been performed starting from a database of about 130 pulses, well representative of the JET discharges with different shape and triangularity of the plasma boundary and with global parameter varying in the whole JET interval. On top of this the base contains several behaviours of the plasma current density profile, starting from monotonic q profile, then reversed shear profile and finally the extreme case of current hole. For some pulses clear MHD signature have been put in evidence and help in particular at the validation of the current density profile or the safety factor q.

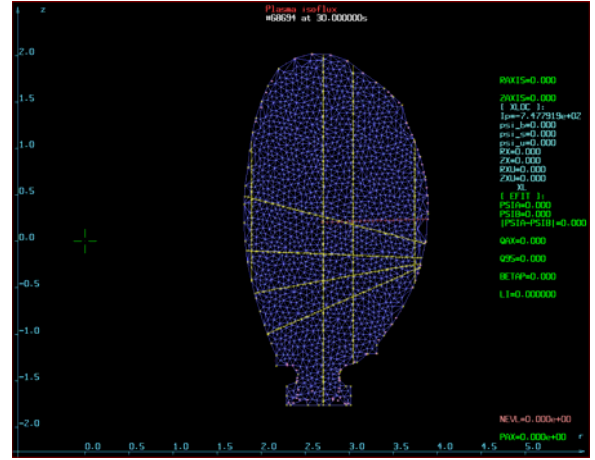

Fig. 2. Equinox Grid used for solving the Grad Shafranov equation

The strategy of the validation has been applied to the two versions of the codes. The first one called Equinox M is the version using only the inputs from the magnetic measurement via XLOC. This version gives accurate plasma geometry and global parameters and does not intend to give very precise information about the q profile. The second version Equinox-J includes internal measurements like polarimetry or MSE and is able to identify hollow plasma current density profiles (validation not showed in this paper)

The validation of the Equinox_M version has been performed mainly using the results of the well assessed EFIT equilibrium code [12] constrained by magnetic measurements only which is used in a routinely manner at JET for intershot analysis. As our code is a free boundary code which means that no assumption is made on the plasma shape, we are able to compare the shape parameters of our reconstruction with the one obtained by XLOC itself. A direct comparison of the plasma Volume and upper triangularity is given in Figure 3 and 4 and shows a very good agreement in terms of the shape of the plasma. In Fig. 5 and Fig. 6 the coordinates Rx and Zx of the X point are compared. We can note in particular a very good agreement of the Rx and some differences for the Zx. Nevertheless these errors are of same order and the origin of the differences is unclear (error on the measurement, method…) . Global quantities like $\beta_p$ and $l_i$ are compared in Fig7. For $\beta_p$ the agreement is quite good and some differences can be noted for the internal inductance. In figure 8 the time evolution of the $\beta_P+l_i/2$ is given for both EFIT and Equinox code and shows again a fairly good agreement. This proves a very good reconstruction of the pressure term and looking at the difference observed for $l_i$ (see fig 9) the medium value is about 0.1. In order to quantify the error and the sensivity of the output related to the error on the measurements we perturbed by 1% the input data of XLOC and get for example a deviation of about 0.1 for the $l_i$. So the difference observed between EFIT and Equinox for the $l_i$ is of the order of the error bars on the results. In terms of q profile this difference is small on the particular case of shot #74937 in Fig10. We can note in particular that the main difference comes from the $q_{ax}$ which is one of the consequences of the lack of information coming from the internal part of the plasma. A very interesting test was then to reconstruct from the obtained equilibrium the line integrated density. In fact once the geometrical line of sight are known and when the plasma equilibrium is obtained, it is possible to solve independently the chi-square minimization of the line integrated density. Results are given in Fig 11 for shot #68690 where a direct comparison between the measured chord 3 and the computed one is given. This agreement is almost perfect.

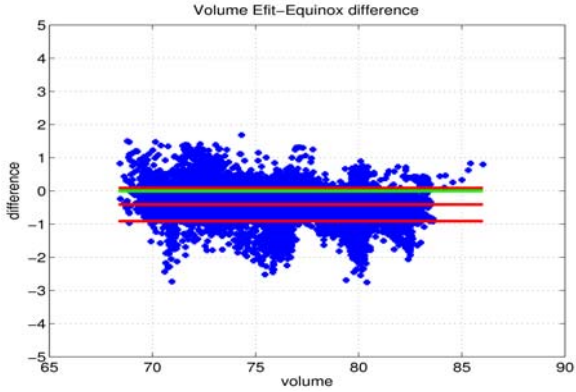

Fig. 3. Standard deviation between EFIT and Equinox for the plasma volume (m3)

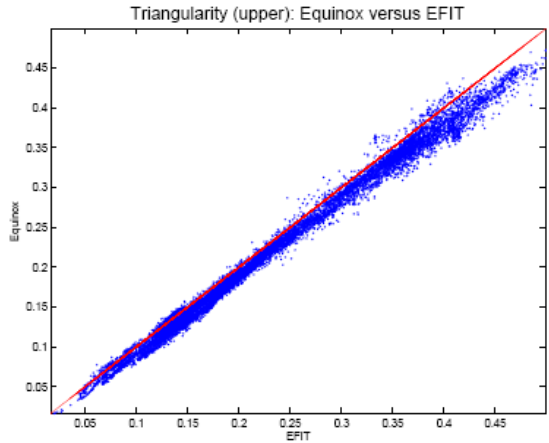

Fig. 4 Comparison between EFIT and Equinox for the upper triangularity

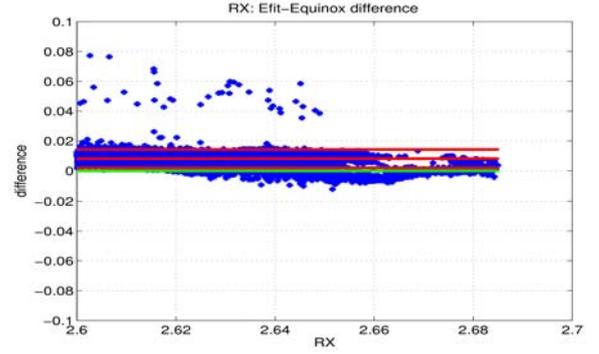

Fig. 5 Comparison between EFIT and Equinox for the R position of the X point

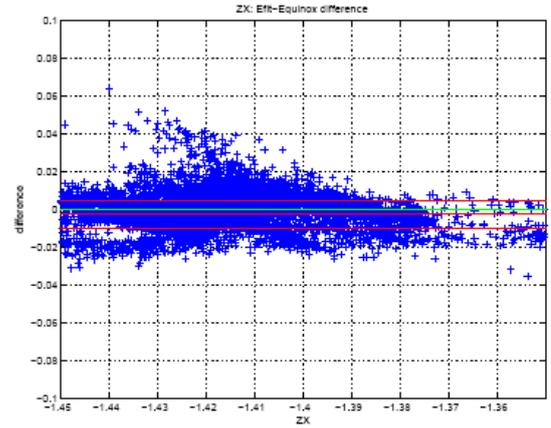

Fig. 6 Comparison between EFIT and Equinox for the Z position of the X point

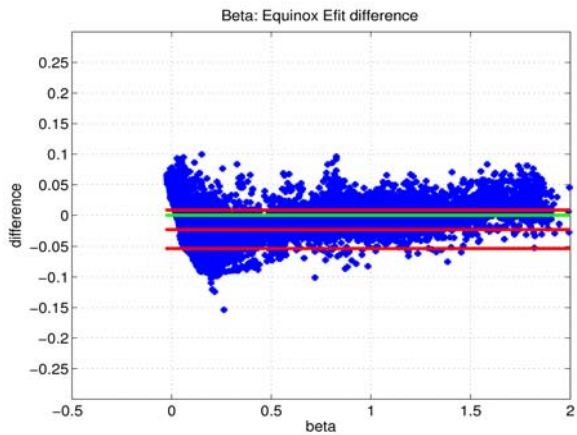

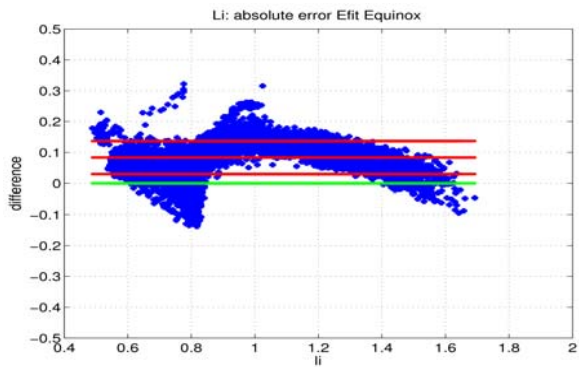

Fig. 7 Comparison between EFIT and Equinox for the beta poloidal (upper figure) and the internal inductance li (lower figure)

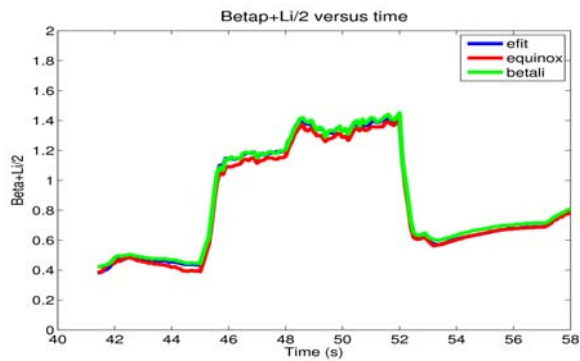

Fig. 8 Comparison between EFIT, Equinox and XLOC of the Shafranov shift.

Finally in order to assess the Equinox reconstruction we have used PROTEUS [13] that solves the direct problem of the Grad shafranov equation. In that particular case a monotonic current density profile was chosen, the equilibrium has been reconstructed and PROTEUS computed the boundary conditions requested by Equinox. Equinox outputs are then compared with the one of PROTEUS. A very good agreement is found on the plasma volume, li and q profile confirming that the statistic relies on a very strong and accurate computation.

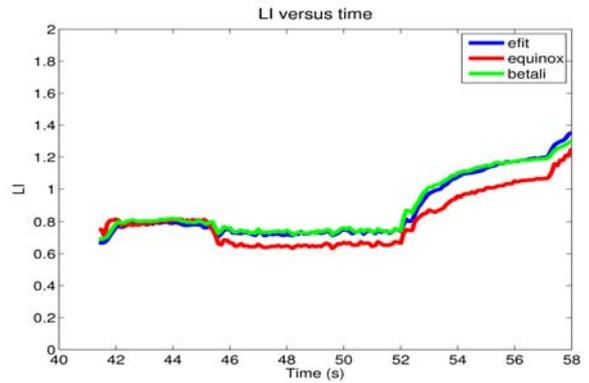

Fig. 9 Comparison between EFIT and Equinox and XLOC of the internal inductance

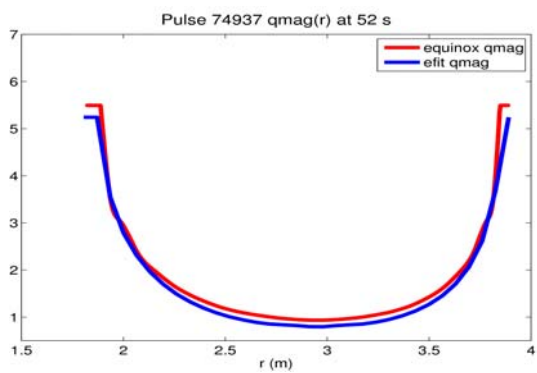

Fig. 10 Comparison between EFIT and Equinox of the safety factor profile

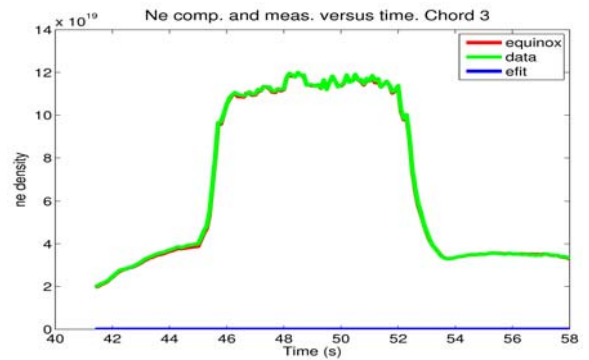

Fig. 11 Comparison between EFIT and Equinox of the line integrated density profile chord3.